\documentclass{amsart}
\usepackage{epsfig}
\usepackage{amssymb}
\newenvironment{example}{\textbf{Example.}}{\par \medskip}
\newenvironment{case}[1]{\nin$\bullet$ \textbf{Case #1.}}{\par}
\newtheorem{theorem}{Theorem}[section]
\newtheorem{prop}[theorem]{Proposition}

\newtheorem{corollary}[theorem]{Corollary}
\newcommand\beq{\begin{equation}}
\newcommand\eeq{\end{equation}}
\newcommand\bce{\begin{center}}
\newcommand\ece{\end{center}}
\newcommand\bea{\begin{eqnarray}}
\newcommand\eea{\end{eqnarray}}
\newcommand\ben{\begin{enumerate}}
\newcommand\een{\end{enumerate}}
\newcommand\brr{\begin{array}}
\newcommand\err{\end{array}}
\newcommand\bt{\begin{tabular}}
\newcommand\et{\end{tabular}}
\newcommand\nin{\noindent}
\newcommand\nn{\nonumber}
\newcommand\bs{\bigskip}
\newcommand\ms{\medskip}

\renewcommand\S{{\mathcal S}}
\def\mn{\mbox{-}}
\def\reduction{\rho}
\def\C{{\mathbf C}}
\def\Ce{{\mathbf C}^{\exp}}
\def\B{{\mathbf B}}
\def\fig{\centering\epsfig}

\title{Asymptotic enumeration of permutations avoiding generalized patterns}
\author{Sergi Elizalde}
\address{Mathematical Sciences Research Institute, 17 Gauss Way, Berkeley, CA 94720}\email{elizalde@msri.org}

\begin{document}
\maketitle \vspace{-12mm}
\begin{abstract}
Motivated by the recent proof of the Stanley-Wilf conjecture, we
study the asymptotic behavior of the number of permutations avoiding
a generalized pattern. Generalized patterns allow the requirement
that some pairs of letters must be adjacent in an occurrence of the
pattern in the permutation, and consecutive patterns are a
particular case of them.

We determine the asymptotic behavior of the number of permutations
avoiding a consecutive pattern, showing that they are an
exponentially small proportion of the total number of permutations.
For some other generalized patterns we give partial results, showing
that the number of permutations avoiding them grows faster than for
classical patterns but more slowly than for consecutive patterns.
\end{abstract}

\section{Introduction}

One of the most important breakthroughs in the subject of
pattern-avoiding permutations has been the proof by Marcus and
Tardos~\cite{MT} of the so-called Stanley-Wilf conjecture, which had
been open for over a decade. This is a basic result regarding the
asymptotic behavior of the number of permutations that avoid a given
pattern. It states that for any pattern $\sigma$ there exists a
constant $\lambda$ such that, if $\alpha_n(\sigma)$ denotes the
number of $\sigma$-avoiding permutations of size $n$, then
$\alpha_n(\sigma)<\lambda^n$. The notion of pattern avoidance that
this result is concerned with is the standard one, namely, where a
permutation is said to avoid a pattern if it does not contain any
subsequence which is order-isomorphic to it.

In \cite{BS}, Babson and Steingr\'{\i}msson introduced the notion of
{\em generalized patterns}, which allows the requirement that
certain pairs of letters of the pattern must be adjacent in any
occurrence of it in the permutation. One particular case of these
are \emph{consecutive patterns}, which were independently studied by
Elizalde and Noy~\cite{EliNoy}. For a subsequence of a permutation
to be an occurrence of a consecutive pattern, its elements have to
appear in adjacent positions of the permutation.

Analogously to the case of classical patterns, it is natural to
study the asymptotic behavior of the number of permutations avoiding
a generalized pattern. This problem is far from being understood. It
follows from our work that for most generalized patterns the number
of permutations avoiding them behaves very differently than in the
case of classical patterns. In this paper we determine the
asymptotic behavior for the case of consecutive patterns, showing
that if $\sigma$ is a consecutive pattern and $\alpha_n(\sigma)$
denotes the number of permutations of size $n$ avoiding it, then
$\lim_{n\rightarrow\infty}\sqrt[n]{\alpha_n(\sigma)/n!}$ is a
positive constant. For some particular generalized patterns we
obtain the same asymptotic behavior, and for patterns of length 3
the problem is solved as well. However, the general case remains
open, and it seems from our investigation that there is a big range
of possible asymptotic behaviors. For some generalized patterns
$\sigma$ of length 4 we give asymptotic upper and lower bounds on
$\alpha_n(\sigma)$.

The paper is structured as follows. In Section~\ref{sec:prel} we
introduce the definitions and notation for generalized pattern
avoidance. We also mention some generating function techniques that
will be used in the paper, as well as previous results regarding
consecutive patterns. In Section~\ref{sec:1sigma} we give the
exponential generating functions for permutations avoiding a special
kind of generalized patterns, extending the results
from~\cite{EliNoy}. In Section~\ref{sec:asymp} we study the
asymptotic behavior as $n$ goes to infinity of the number of
permutations of size $n$ avoiding a generalized pattern, solving the
problem only in some cases. In Section~\ref{sec:12_34} we give lower
and upper bounds on the number of $12\mn34$-avoiding permutations,
and in Section~\ref{sec:1_23_4} we obtain a similar result for the
pattern $1\mn23\mn4$. Finally, in Section~\ref{sec:other} we discuss
some open problems and further research.

\section{Preliminaries}\label{sec:prel}

In this section we define most of the notation that will be used
later on. We start introducing the notion of generalized pattern
avoidance.

\subsection{Generalized patterns}\label{sec:gp}
These patterns, which were introduced by Babson and
Steingr\'{\i}msson~\cite{BS}, extend the classical notion of pattern
avoidance. We will denote by $\S_n$ the symmetric group on
$\{1,2,\ldots,n\}$. Let $n$, $m$ be two positive integers with $m\le
n$, and let $\pi=\pi_1\pi_2\cdots\pi_n\in\S_n$ be a permutation. A
generalized pattern $\sigma$ is obtained from a permutation
$\sigma_1\sigma_2\cdots\sigma_m\in\S_m$ by choosing, for each
$j=1,\ldots,m-1$, either to insert a dash $\mn$ between $\sigma_j$
and $\sigma_{j+1}$ or not. More formally,
$\sigma=\sigma_1\varepsilon_1\sigma_2\varepsilon_2\cdots\varepsilon_{m-1}\sigma_m$,
where each $\varepsilon_j$ is either the symbol~$\mn$ or the empty
string.

With this notation, we say that $\pi$ \emph{contains} (the
generalized pattern) $\sigma$ if there exist indices
$i_1<i_2<\ldots<i_m$ such that
\renewcommand{\theenumi}{\roman{enumi}} \ben
\item for each $j=1,\ldots,m-1$, if $\varepsilon_j$ is empty
then $i_{j+1}=i_j+1$, and
\item $\reduction(\pi_{i_1}\pi_{i_2}\cdots\pi_{i_m})=\sigma_1\sigma_2\cdots\sigma_m$,
where $\reduction$ is the reduction consisting in relabeling the
elements with $\{1,\ldots,m\}$ so that they keep the same order
relationships they had in $\pi$. (Equivalently, this means that for
all indices $a$ and $b$, $\pi_{i_a}<\pi_{i_b}$ if and only if
$\sigma_a<\sigma_b$.) \een In this case,
$\pi_{i_1}\pi_{i_2}\cdots\pi_{i_m}$ is called an \emph{occurrence}
of $\sigma$ in $\pi$.

If $\pi$ does not contain $\sigma$, we say that $\pi$ \emph{avoids}
$\sigma$, or that it is \emph{$\sigma$-avoiding}. For example, the
permutation $\pi=3542716$ contains the pattern $12\mn4\mn3$, and it
has exactly one occurrence of it, namely the subsequence $3576$. On
the other hand, $\pi$ avoids the pattern $12\mn43$.

Observe that in the case where $\sigma$ has dashes in all $m-1$
positions, we recover the classical definition of pattern avoidance,
because in this case condition (i) holds trivially. On the other
end, the case in which $\sigma$ has no dashes corresponds to
consecutive patterns. In this situation, an occurrence of $\sigma$
in $\pi$ has to be a consecutive subsequence. Consecutive patterns
were introduced independently in \cite{EliNoy}, where the authors
give generating functions for the number of occurrences of certain
consecutive patterns in permutations. Several papers deal with the
enumeration of permutations avoiding generalized patterns. In
\cite{C}, Claesson presented a complete solution for the number of
permutations avoiding any single $3$-letter generalized pattern with
exactly one adjacent pair of letters. Claesson and Mansour \cite{CM}
(see also \cite{Mg1}) did the same for any pair of such patterns.
In~\cite{EliMan}, Elizalde and Mansour studied the distribution of
several statistics on permutations avoiding $1\mn3\mn2$ and $1\mn23$
simultaneously. On the other hand, Kitaev \cite{Ki} investigated
simultaneous avoidance of two or more $3$-letter generalized
patterns without dashes.

All the patterns that appear in this paper will be represented by
the notation just described. In particular, a pattern
$\sigma=\sigma_1\sigma_2\cdots\sigma_m$ without dashes will denote a
consecutive pattern. We will represent classical patterns by writing
dashes between any two adjacent elements, namely, as
$\sigma_1\mn\sigma_2\mn\cdots\mn\sigma_m$.

If $\sigma$ is a generalized pattern, let $\S_n(\sigma)$ denote the
set of permutations in $\S_n$ that avoid $\sigma$. Let
$\alpha_n(\sigma)=|\S_n(\sigma)|$ be the number of such
permutations, and let
$$A_\sigma(z)=\sum_{n\ge0}\alpha_n(\sigma)\frac{z^n}{n!}$$ be the
exponential generating function counting $\sigma$-avoiding
permutations.

\subsection{Labeled classes and exponential generating
functions}\label{sec:labeled}

Here we recall some basic machinery for exponential generating
functions that will be used later. We direct the reader to
\cite{FlSe98} for a detailed account on combinatorial classes and
the symbolic method. Let ${\mathcal A}$ be a class of labeled
combinatorial objects and let $|\zeta|$ be the size of an object
$\zeta \in {\mathcal A}$. If ${\mathcal A}_n$ denotes the objects in
${\mathcal A}$ of size $n$ and $a_n = |{\mathcal A}_n|$, then the
{\em exponential generating function\/}, EGF for short, of the class
${\mathcal A}$ is
 $$
     A(z) = \sum_{\zeta\in {\mathcal A}} \frac{z^{|\zeta|}}{|\zeta|!}
     = \sum_{n\ge0} a_n \frac{z^n}{n!}.
 $$
In our context, the size of a permutation is simply its length.

There is a direct correspondence between set-theoretic operations
(or ``constructions'') on combinatorial classes and algebraic
operations on EGFs. Table~\ref{tab:symb} summarizes this
correspondence for the operations that are used in the paper. There
``union'' means union of disjoint copies, ``labeled product'' is the
usual cartesian product enriched with the relabeling operation, and
``set'' forms sets in the usual sense. Particularly important for us
is the construction ``boxed product'' ${\mathcal A}={\mathcal
B}^\Box*{\mathcal C}$, which corresponds to the subset of
$\mathcal{B}\star\mathcal{C}$ (the usual labeled product) formed by
those pairs in which the smallest label lies in the $\mathcal{B}$
component. Another similar construction is the ``double boxed
product'' ${\mathcal A}={\mathcal B}^\boxtimes*{\mathcal C}$, which
denotes the subset of $\mathcal{B}\star\mathcal{C}$ formed by those
pairs in which both the smallest and the largest label lie in the
$\mathcal{B}$ component.

\renewcommand{\arraystretch}{1.3}
\begin{table}[htb]
\centering
\begin{tabular}{l|l|l} \hline
\em Construction &  &  \em Operation on GF \\
 \hline
 Union & ${\mathcal A}={\mathcal B} \cup {\mathcal C}$  &  $A(z)=B(z)+C(z)$ \\
 Labeled product & ${\mathcal A}={\mathcal B}\star{\mathcal C}$  & $A(z)=B(z)C(z)$ \\
 Set & ${\mathcal A} = \Pi ({\mathcal B})$ &  $A(z) = \exp(B(z))$ \\
 Boxed product &  ${\mathcal A}={\mathcal B}^\Box \star {\mathcal C}$
    & $A(z)=\int_0^{z}(\frac{d}{dt}B(t))\cdot C(t)\, dt$\\
 Double boxed product &  ${\mathcal A}={\mathcal B}^\boxtimes \star {\mathcal C}$
    & $A(z)=\int_0^{z}\int_0^{y}(\frac{d^2}{dt^2}B(t))\cdot C(t)\, dt\,dy$\\
 \hline
\end{tabular}
\ms \caption{\label{tab:symb} The basic combinatorial constructions
and their translation into exponential generating functions.}
\end{table}

\subsection{Consecutive patterns of length $3$}

For patterns of length $3$ with no dashes, it follows from the
trivial reversal and complementation operations that
$\alpha_n(123)=\alpha_n(321)$ and
$\alpha_n(132)=\alpha_n(231)=\alpha_n(312)=\alpha_n(213)$. The EGFs
for these numbers are given in the following theorem of Elizalde and
Noy~\cite{EliNoy}, which we will use later in the paper. The symbol
$\sim$ between two sequences indicates that they have the same
asymptotic behavior.

\begin{theorem}[\cite{EliNoy}]\label{th:EN}
We have $$A_{123}(z)
=\frac{\sqrt{3}}{2}\frac{e^{z/2}}{\cos(\frac{\sqrt{3}}{2}z+
\frac{\pi}{6})}, \hspace{2cm} A_{132}(z)=\frac{1}{1-\int_0^z
e^{-t^2/2} dt}. $$
 Their coefficients satisfy
$$\alpha_n(123)\sim\gamma_1\cdot(\rho_1)^n\cdot n!, \hspace{2cm}
\alpha_n(132)\sim\gamma_2\cdot(\rho_2)^{n}\cdot n!, $$
 where $\rho_1=\frac{3\sqrt{3}}{2\pi}$,
$\gamma_1=e^{3\sqrt{3}\pi}$, $(\rho_2)^{-1}$ is the unique positive
root of $\int_0^z e^{-t^2/2} dt=1$, and
$\gamma_2=e^{({\rho_2})^{-2}/2}$, the approximate values being
$$\rho_1=0.8269933, \quad\gamma_1=1.8305194, \quad\rho_2=0.7839769,
\quad\gamma_2=2.2558142.$$ Furthermore, for every $n\geq 4$, we have
$$\alpha_n(123)>\alpha_n(132).$$
\end{theorem}

\section{Patterns of the form $1\mn\sigma$}\label{sec:1sigma}

In this section we study a very particular class of generalized
patterns, namely those that start with $1\mn$, followed by a
consecutive pattern (i.e., without dashes).

\begin{prop}\label{prop:1sigma}
Let $\sigma=\sigma_1\sigma_2\cdots\sigma_k\in\S_k$ be a
consecutive pattern, and let $1\mn\sigma$ denote the generalized
pattern $1\mn(\sigma_1+1)(\sigma_2+1)\cdots(\sigma_k+1)$. Then,
$$A_{1\mn\sigma}(z)=\exp\left(\int_0^z A_\sigma(t)\, dt\right).$$
\end{prop}

\begin{proof}
Given a permutation $\pi$, let $m_1>m_2>\cdots>m_r$ be the values of
its left-to-right minima (recall that $\pi_i$ is a left-to-right
minimum of $\pi$ if $\pi_j>\pi_i$ for all $j<i$). We can write
$\pi=m_1 w_1 m_2 w_2 \cdots m_r w_r$, where each $w_i$ is a
(possibly empty) subword of $\pi$, each of whose elements is greater
than $m_i$. We claim that $\pi$ avoids $1\mn\sigma$ if and only if
each of the blocks $w_i$ (more precisely, its reduction
$\reduction(w_i)$) avoids the consecutive pattern $\sigma$. Indeed,
it is clear that if one of the blocks $w_i$ contains $\sigma$, then
$m_i$ together with the occurrence of $\sigma$ forms an occurrence
of $1\mn\sigma$. Conversely, if $\pi$ contains $1\mn\sigma$, then
the elements of $\pi$ corresponding to $\sigma$ have to be adjacent,
and none of them can be a left-to-right minimum (since the element
corresponding to `1' has to be to their left), therefore they must
be all inside the same block $w_i$ for some $i$.

If we denote by ${\mathcal A}$ the class of permutations avoiding
$\sigma$, then, in the notation of Table~\ref{tab:symb}, the class
of permutations avoiding $1\mn\sigma$ can be expressed as
$$\Pi (\{ z\}^\Box \star {\mathcal A}),$$ where $\{ z\}^\Box \star
{\mathcal A}$ corresponds to a block $m_i w_i$, with the box
indicating that the left-to-right minimum has the smallest label.
The set construction arises from the fact given a collection of
blocks $m_i w_i$, there is a unique way to order them, namely with
the left-to-right minima in decreasing order. The expression
$A_{1\mn\sigma}(z)=\exp(\int_0^z A_\sigma(t)\, dt)$ follows now from
this construction.
\end{proof}

Proposition~\ref{prop:1sigma} also appears independently in a
preprint of Kitaev~\cite{Ki2}.\ms

\begin{example}
The only permutation avoiding $\sigma=12$ (resp. $\sigma=21$) is the
decreasing (resp. increasing) one. Therefore, by
Proposition~\ref{prop:1sigma},
$$A_{1\mn23}(z)=A_{1\mn32}(z)=\exp\left(\int_0^z e^t dt\right)=e^{e^z-1},$$
the EGF for Bell numbers, which agrees with the result
in~\cite{C}.
\end{example}

\begin{example}
For the consecutive patterns $132$, $231$, $312$ and $213$, the
generating function for the number of permutations avoiding either
of them is given in Theorem~\ref{th:EN} (which follows from
\cite[Theorem 4.1]{EliNoy}). Now, by Proposition~\ref{prop:1sigma},
we get the following expression:
$$A_{1\mn243}(z)=A_{1\mn342}(z)=A_{1\mn423}(z)=A_{1\mn324}(z)=\exp\left(\int_0^z\frac{dt}{1-\int_0^z e^{-u^2/2} du}\right).$$
\end{example}

\begin{example}
The EGF for permutations avoiding $123$ or $321$ is also given in
Theorem~\ref{th:EN}. Proposition~\ref{prop:1sigma} implies now that
$$A_{1\mn234}(z)=A_{1\mn432}(z)=\exp\left(\frac{\sqrt{3}}{2}\int_0^z
\frac{e^{t/2}\,dt}{\cos(\frac{\sqrt{3}}{2}t+\frac{\pi}{6})}\right).$$
\end{example}

Combined with the results of \cite{EliNoy},
Proposition~\ref{prop:1sigma} gives expressions for the EGFs
$A_{1\mn\sigma}(z)$ where $\sigma$ has one of the following forms:

\bt{l} $\sigma=123\cdots k$,\\
$\sigma=k(k-1)\cdots21$,\\ $\sigma=12\cdots a\,\tau\,(a+1)$,\\
$\sigma=(a+1)\,\tau\,a(a-1)\cdots21$,\\
$\sigma=k(k-1)\cdots(k+1-a)\,\tau'\,(k-a)$,\\
$\sigma=(k-a)\,\tau'\,(k+1-a)(k+2-a)\cdots k$, \et

where $k,a$ are positive integers with $a\le k-2$, $\tau$ is any
permutation of $\{a+2,a+3,\cdots,k\}$ and $\tau'$ is any
permutation of $\{1,2,\cdots,k-a-1\}$.

\section{Asymptotic enumeration}\label{sec:asymp}

Here we discuss the behavior of the numbers $\alpha_n(\sigma)$ as
$n$ goes to infinity, for a given generalized pattern $\sigma$. We
use the symbol $\sim$ to indicate that two sequences of numbers
have the same asymptotic behavior (i.e., we write $a_n\sim b_n$ if
$\lim_{n\rightarrow\infty}\frac{a_n}{b_n}=1$), and we use the
symbol $\ll$ to indicate that a sequence is asymptotically smaller
than another one (i.e., we write $a_n\ll b_n$ if
$\lim_{n\rightarrow\infty}\frac{a_n}{b_n}=0$).

Let us first consider the case of consecutive patterns.

\begin{theorem}\label{asympconsec}
Let $k\ge3$ and let $\sigma\in\S_k$ be a consecutive pattern. \ben
\item There exist constants $0<c,d<1$ such that $$c^n n! <
\alpha_n(\sigma) < d^n n!$$ for all $n\ge k$.
\item There exists a constant $0< w \le 1$ such that
$$\lim_{n\rightarrow\infty}\left(\frac{\alpha_n(\sigma)}{n!}\right)^{1/n} = w.$$
\een
\end{theorem}

Note that $c$, $d$ and $w$ depend only on $\sigma$. Compare this
result with the conjecture of Warlimont~\cite{War} that for any
consecutive pattern $\sigma$ there exist constants $\gamma>0$ and
$w<1$ such that $\alpha_n(\sigma)/n!\sim \gamma\, w^n$.

\begin{proof}
The key observation is that, for any consecutive pattern $\sigma$,
\beq\label{eq:mn} \alpha_{m+n}(\sigma) \le
\alpha_{m}(\sigma)\alpha_{n}(\sigma)
     \binom{m+n}{n}.\eeq
To see this, just observe that a $\sigma$-avoiding permutation of
length $m+n$ induces two juxtaposed $\sigma$-avoiding permutations
of lengths $m$ and $n$.

By induction on $n \ge k$ one gets
 $$
    \alpha_{m+n}(\sigma) < d^m m! \, d^n n! \binom{m+n}{n} =
    d^{m+n}(m+n)!
 $$
for some positive $d<1$.

For the lower bound, let $\tau=\reduction(\sigma_1\sigma_2\sigma_3)$
be the reduction of the first three elements of $\sigma$. Clearly
$\S_n(\tau)\subseteq\S_n(\sigma)$ for all $n$, since an occurrence
of $\sigma$ in a permutation produces also an occurrence of $\tau$,
hence $\alpha_n(\tau)\leq\alpha_n(\sigma)$. But the fact that
$\sigma\in\S_3$ implies that $\alpha_n(\sigma)$ equals either
$\alpha_n(123)$ or $\alpha_n(132)$. In any case, by
Theorem~\ref{th:EN} we have that
$$\alpha_n(\sigma)\geq\alpha_n(132)> c^n n!$$ for some $c>0$.

To prove part (ii), we can express~(\ref{eq:mn}) as
$$\frac{\alpha_{m+n}(\sigma)}{(m+n)!} \le
\frac{\alpha_{m}(\sigma)}{m!}\frac{\alpha_{n}(\sigma)}{n!}$$ and
apply {\it Fekete's lemma} (see \cite[Lemma 11.6]{vLW}
or~\cite{Fek}) to the function $n!/\alpha_{n}(\sigma)$ to conclude
that
$\lim_{n\rightarrow\infty}\left(\frac{\alpha_n(\sigma)}{n!}\right)^{1/n}$
exists. Calling it $w$, then part~(i) implies that $w\le 1$ and
$w\ge\lim_{n\rightarrow\infty}\left(\frac{\alpha_n(132)}{n!}\right)^{1/n}=0.7839769$.
\end{proof}

\ms

In order to study the asymptotic behavior of $\alpha_n(\sigma)$ for
a generalized pattern $\sigma$ we separate the problem into the
following three cases. Assume from now on that $k\ge3$ and that
$\sigma$ is a generalized pattern of length $k$. We use the word
\emph{slot} to refer to the place between two adjacent elements of
$\sigma$, where there can be a dash or not.

\bs

\begin{case}{1}
{\it The pattern $\sigma$ has dashes between any two adjacent
elements, i.e., $\sigma=\sigma_1\mn\sigma_2\mn\cdots\mn\sigma_k$.}
\ms

These are just the classical patterns, which have been widely
studied in the literature. The asymptotic behavior of the number of
permutations avoiding them is given by the Stanley-Wilf conjecture,
which has been recently proved by Marcus and Tardos~\cite{MT}, after
several authors had given partial results over the last few
years~\cite{AF,Arr,B2,Kla}.

\begin{theorem}[\emph{Stanley-Wilf conjecture}, proved
in~\cite{MT}] \label{th:swconj} For every classical pattern
$\sigma=\sigma_1\mn\sigma_2\mn\cdots\mn\sigma_k$, there is a
constant $\lambda$ (which depends only on $\sigma$) such that
$$\alpha_n(\sigma)<\lambda^n$$ for all $n\ge1$.
\end{theorem}

On the other hand, it is clear that
$\alpha_n(\sigma)\ge\alpha_n(\reduction(\sigma_1\mn\sigma_2\mn\sigma_3))=\C_n\sim\frac{1}{\sqrt{\pi
n}}\,4^n$, where $\C_n$ denotes the $n$-th \emph{Catalan number}. As
shown by Arratia~\cite{Arr}, Theorem~\ref{th:swconj} is equivalent
to the statement that
$\lim_{n\rightarrow\infty}\sqrt[n]{\alpha_n(\sigma)}$ exists. The
value of this limit has been computed for several classical
patterns: it is clearly $4$ for patterns of length $3$, it is
known~\cite{Reg} to be $(k-1)^2$ for $\sigma=1\mn2\mn\cdots\mn k$,
it has been shown~\cite{Bon97} to be $8$ for $\sigma=1\mn3\mn4\mn2$,
and it has recently been proved by B\'ona~\cite{B3} to be
nonrational for certain patterns.
\end{case}

\ms

\begin{case}{2}
{\it The pattern $\sigma$ has two consecutive slots without a dash
(equivalently, three consecutive elements without a dash between
them), i.e., $\sigma=\cdots\sigma_i\sigma_{i+1}\sigma_{i+2}\cdots$.}
\ms

\begin{prop} Let $\sigma$ be a generalized pattern having three consecutive elements without a
dash. Then there exist constants $0<c,d<1$ such that $$c^n n! <
\alpha_n(\sigma) < d^n n!$$ for all $n\ge k$.
\end{prop}

\begin{proof} For the upper bound, notice that if a permutation
contains the consecutive pattern
$\sigma_1\sigma_2\sigma_3\cdots\sigma_k$ obtained by removing all
the dashes in $\sigma$, then it also contains $\sigma$. Therefore,
$\alpha_n(\sigma)\le\alpha_n(\sigma_1\sigma_2\sigma_3\cdots\sigma_k)$
for all $n$, and now the upper bound follows from part~(i) of
Theorem~\ref{asympconsec}.

For the lower bound, we use that
$\alpha_n(\sigma)\ge\alpha_n(\reduction(\sigma_i\sigma_{i+1}\sigma_{i+2}))
\ge\alpha_n(132)> c^n n!$, where $\sigma_i\sigma_{i+1}\sigma_{i+2}$
are three consecutive elements in $\sigma$ without a dash.
\end{proof}
\end{case}

\ms

\begin{case}{3}
{\it The pattern $\sigma$ has at least a slot without a dash, but
not two consecutive slots without dashes.} \ms

This case includes all the patterns not considered in Cases 1 and~2.
The asymptotic behavior of $\alpha_n(\sigma)$ for these patterns is
not known in general. The case of patterns of length~3 is covered by
the following result, due to Claesson~\cite{C}. Let $\B_n$ denote
the $n$-th \emph{Bell number}, which counts the number of partitions
of an $n$-element set.

\begin{prop}[\cite{C}]\label{prop:claesson} Let $\sigma$ be a generalized pattern of length~3 with one dash. \ben
\item If
$\sigma\in\{1\mn23,3\mn21,32\mn1,12\mn3,1\mn32,23\mn1,3\mn12,21\mn3\}$,
then $\alpha_n(\sigma)=\B_n$.
\item If $\sigma\in\{2\mn13,2\mn31,31\mn2,13\mn2\}$, then $\alpha_n(\sigma)=\C_n$.
\een
\end{prop}

It is known that the asymptotic behavior of the Catalan numbers is
given by $\C_n\sim\frac{1}{\sqrt{\pi n}}\,4^n$. For the Bell
numbers, one has the formula
$$\B_n\sim\frac{1}{\sqrt{n}}\,\lambda(n)^{n+1/2}e^{\lambda(n)-n-1},$$
where $\lambda(n)$ is defined by $\lambda(n)\ln(\lambda(n))=n$.
Another useful description of the asymptotic behavior of $\B_n$ is
the following: $$\frac{\ln\B_n}{n}=\ln n-\ln\ln
n+O\left(\frac{\ln\ln n}{\ln n}\right).$$ This shows in particular
that $\delta^n\ll\B_n\ll c^n n!$ for any constants $\delta,c>0$.

For patterns $\sigma$ of length $k\ge4$ that have slots without a
dash, but not two consecutive slots without dashes, not much is
known in general about the number of permutations avoiding them. It
follows from Cases 1 and 2 that
$$\delta^n<\alpha(\sigma_1\mn\sigma_2\mn\cdots\mn\sigma_k)\le\alpha_n(\sigma)
\le\alpha(\sigma_1\sigma_2\cdots\sigma_k)<d^n n!$$ for some
constants $\delta>0$ and $d<1$. Clearly, if $\sigma$ contains one of
the patterns in part~(i) of Proposition~\ref{prop:claesson}, then
the lower bound can be improved to $\B_n$. However, determining the
precise asymptotic behavior of $\alpha_n(\sigma)$ seems to be a
difficult problem. In the rest of the paper we discuss a few partial
results in this direction.
\end{case}

\ms

The next statement is about permutations of the form $1\mn\sigma$.

\begin{corollary} Let $\sigma$ be a consecutive pattern, and let
$1\mn\sigma$ be defined as in Proposition~\ref{prop:1sigma}. Then,
$$\lim_{n\rightarrow\infty}\left(\frac{\alpha_n(1\mn\sigma)}{n!}\right)^{1/n}
=\lim_{n\rightarrow\infty}\left(\frac{\alpha_n(\sigma)}{n!}\right)^{1/n}.$$
\end{corollary}

\begin{proof}
By Proposition~\ref{prop:1sigma} we know that
$A_{1\mn\sigma}(z)=\exp\left(\int_0^z A_\sigma(t)\, dt\right)$.
Since the exponential is an analytic function over the whole complex
plane, we obtain that $A_{1\mn\sigma}(z)$ has the same radius of
convergence as $A_{\sigma}(z)$, from where the result follows.
\end{proof}

\section{The pattern $12\mn34$}\label{sec:12_34}

The next proposition gives an upper and a lower bound for the
numbers $\alpha_n(12\mn34)$. Given two formal power series
$F(z)=\sum_{n\ge0}f_nz^n$ and $G(z)=\sum_{n\ge0}g_nz^n$, we use the
notation $F(z)< G(z)$ to indicate that $f_n<g_n$ for all $n$, and
$F(z)\ll G(z)$ to indicate that $f_n\ll g_n$.

\begin{prop}\label{prop:12_34}
For $k\ge1$, let \bea\nn \hspace*{-1cm} && h_k=1+\frac{1}{2}+\cdots+\frac{1}{k}, \\
\nn && b_k(z)=\sum_{i=0}^k\binom{k}{i}^2[z+2(h_{k-i}-h_i)]\,e^{iz}, \\
\nn &&
c_k(z)=\frac{e^{(k+1)z}}{k+1}-\sum_{i=0}^k\binom{k}{i}\binom{k+1}{i}
\left[z+2(h_{k-i}-h_i)+\frac{1}{k+1-i}\right]e^{iz}, \\
\nn && S(z)=\sum_{k\ge1}b_k(z)+\sum_{k\ge1}c_k(z). \eea Then
$$e^{S(z)}< A_{12\mn34}(z)< e^{S(z)+e^z+z-1}.$$
\end{prop}

\begin{figure}[hbt] \fig{file=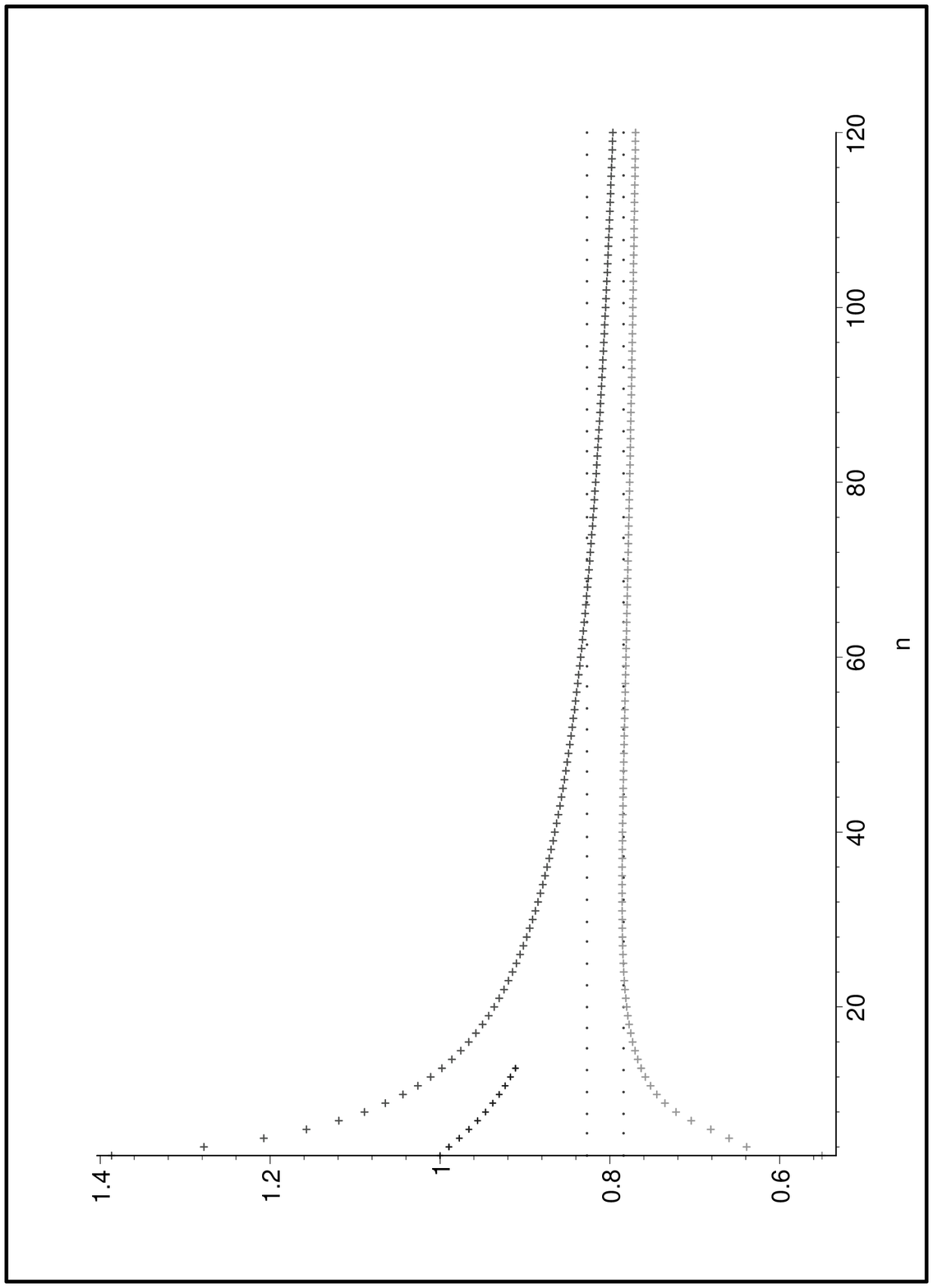,height=12.5cm,angle=-90}
\caption{\label{fig:bounds} The first values of
$\sqrt[n]{\alpha_n(12\mn34)/n!}$ between the lower and the upper
bound given by Proposition~\ref{prop:12_34}.}
\end{figure}

If we write $e^{S(z)}=\sum l_n \frac{z^n}{n!}$ and
$e^{S(z)+e^z+z-1}=\sum u_n \frac{z^n}{n!}$ to denote the
coefficients of the series giving the lower and the upper bound
respectively, then the graph in Figure~\ref{fig:bounds} shows the
values of $\sqrt[n]{\alpha_n(12\mn34)/n!}$ for $n\le13$, bounded
between the values $\sqrt[n]{l_n/n!}$ and $\sqrt[n]{u_n/n!}$ for
$n\le120$. The two horizontal dotted lines are at height $0.7839769$
and $0.8269933$, which are
$\lim_{n\rightarrow\infty}\sqrt[n]{\alpha_n(\sigma)/n!}$ for
$\sigma=132$ and $\sigma=123$ respectively, given by
Theorem~\ref{th:EN}. From this plot it seems conceivable that
$\lim_{n\rightarrow\infty}\sqrt[n]{\alpha_n(12\mn34)/n!}=0$,
although we have not succeeded in proving this.

Note that the lower bound, together with the fact that $S(z)\gg
e^z-1$ (which follows from the definition), shows that
$A_{12\mn34}(z)> e^{S(z)} \gg e^{e^z-1}$, which means that
$\alpha_n(12\mn34)\gg\B_n$, that is, the number of
$12\mn34$-avoiding permutations is asymptotically larger than the
Bell numbers.

\begin{proof}
Let $\pi$ be a permutation that avoids $12\mn34$. This means that it
has no two ascents such that the second one starts at a higher value
than where the first one ends. We can write
$\pi=B_0a_1B_1a_2B_2a_3B_3\cdots$, where $a_1$ and the element
preceding it form the first ascent of $\pi$, $a_2$ and the element
preceding it form the first ascent such that $a_2<a_1$, $a_3$ and
the element preceding it form the first ascent such that $a_3<a_2$,
and so on. By definition, $B_0$ is a non-empty decreasing word whose
last element is less than $a_1$, and each $B_i$ with $i\ge1$ can be
written uniquely as a sequence $B_i=w_{i,0}w_{i,1}w_{i,2}\cdots
w_{i,r_i}$ for some $r_i\ge1$ ($r_i$ can be 0 if $w_{i,0}$ is
nonempty) with the following properties:\ben
\item each $w_{i,j}$ is a decreasing word,
\item for $j\ge1$, $w_{i,j}$ is nonempty and its first element is bigger than
$a_i$,
\item the last element of each $w_{i,j}$ is less than $a_i$,
\item the last element of $B_i$ is less than $a_{i+1}$.
\een These properties ensure that $\pi$ avoids $12\mn34$ (since no
$B_i$ has an ascent above $a_i$), and that the decomposition is
unique.

Ideally we would like to use this decomposition to find a
generating function for the numbers $\alpha_n(12\mn34)$.
Unfortunately, the structure of the decomposition is a bit too
complicated to find an exact formula. Instead, we will add and
remove restrictions to simplify this description, which allows us
to give lower and upper bounds respectively.

\ms

To find an upper bound, we will count permutations of the form
$\pi=B_0a_1B_1a_2B_2a_3B_3\cdots$, where the $B_i$ and $a_i$ satisfy
the properties above, except for the requirement (iv) that the last
element of each $B_i$ has to be less than $a_{i+1}$. Omitting this
requirement we are overcounting permutations, and thus we get an
upper bound. The first step now is to find the EGF for a block $K_i$
of the form $a_iB_i$, where $B_i$ satisfies properties (i), (ii) and
(iii) from above.

Let us first assume that $w_{i,0}$ is empty, that is,
$B_i=w_{i,1}w_{i,2}\cdots w_{i,r_i}$. We compute the EGF for
$K_i=a_iB_i$ where $r_i$ is fixed, by induction on $r_i$. If
$r_i=0$, then we have that $K_i=a_i$, so the EGF is $b_0(z):=z$. If
$r_i=1$, then $K_i=a_iw_{i,1}$, where $w_{i,1}$ is a decreasing word
starting above $a_i$ and ending below it. The EFG for $w_{i,1}$ is
$e^z$. Now, to incorporate the condition that the largest and the
smallest labels of $K_i$ lie in $w_{i,1}$, we use the double boxed
product construction described in Section~\ref{sec:labeled}. A
double derivative is needed to mark the two special elements. We get
that the EGF for such a block is
$$\int_0^z\int_0^y t\left(\frac{d^2}{dt^2}e^t\right) dt\,dy=
\int_0^z\int_0^y t e^t dt\,dy=(z-2)e^z+z+2=b_1(z).$$ Let now
$r_i=2$. The case in which both the largest and the smallest label
of $K_i=a_iw_{i,1}w_{i,2}$ are contained in $w_{i,2}$ corresponds to
the EGF \beq\label{maxminw2} \int_0^z\int_0^y
b_1(t)\left(\frac{d^2}{dt^2}e^t\right) dt\,dy.\eeq If we write each
$w_{i,j}$ as $w^+_{i,j}w^-_{i,j}$, separating the elements above and
below $a_i$ ($w^+_{i,j}$ and $w^-_{i,j}$ respectively), then the
largest element of $K_i$ can be either in $w^+_{i,1}$ or in
$w^+_{i,2}$, and the smallest element of $K_i$ can be either in
$w^-_{i,1}$ or in $w^-_{i,2}$. Thus, all the possibilities are
obtained from the case counted by the EGF (\ref{maxminw2}) by
permuting the upper and lower parts of the $w_{i,1}$ and $w_{i,2}$
in the four possible different ways. It follows that the EGF for
$K_i$ when $r_i=2$ is
$$4\int_0^z\int_0^y b_1(t)\,e^t
dt\,dy=(z-3)e^{2z}+4ze^z+z+3=b_2(z).$$ In general, if $b_{k-1}(z)$
is the EGF for the case $r_i=k-1$, then the EGF for the case
$r_i=k$ is given by $$b_k(z)=k^2\int_0^z\int_0^y b_{k-1}(t)\,e^t
dt\,dy.$$ It is straightforward to check that the functions
$b_k(z)$ defined in the statement of the proposition satisfy this
recurrence.

The case where $w_{i,0}$ is nonempty can be treated similarly. Now
we have $B_i=w_{0,i}w_{i,1}w_{i,2}\cdots w_{i,r_i}$. If $r_i=0$,
the EGF for $a_iw_{0,i}$ is $c_0(z):=e^z-1-z$ (since the block has
at least 2 elements). If $r_i=1$, then a block of the form
$a_iw_{0,i}w_{i,1}$ can be obtained from the case where the
largest and the smallest element are in $w_{i,1}$ by permuting
$w_{0,i}$ and $w^-_{i,1}$ if necessary. This yields the EGF
$$2\int_0^z\int_0^y c_0(t)\left(\frac{d^2}{dt^2}e^t\right)
dt\,dy=\frac{e^{2z}}{2}+2(1-z)e^z-z-\frac{5}{2}=c_1(z).$$ In
general, for nonempty $w_{i,0}$, if $c_{k-1}(z)$ is the EGF for
the case $r_i=k-1$, then the EGF for the case $r_i=k$ is given by
$$c_k(z)=k(k+1)\int_0^z\int_0^y c_{k-1}(t)\,e^t dt\,dy.$$
This is the recurrence satisfied by the functions $c_k(z)$ defined
in the statement of the proposition.

The generating function for a set of blocks $K_i=a_iB_i$ of the
form just described is
$$\exp\left(\sum_{k\ge0}b_k(z)+\sum_{k\ge0}c_k(z)\right)=\exp(S(z)+z+e^z-1-z).$$
From such a set there is a unique way to form a sequence
$a_1B_1a_2B_2a_3B_3\cdots$ where $a_1>a_2>a_3>\cdots$. Finally, we
multiply by $e^z$ to take into account the initial decreasing
segment $B_0$ of the permutation $\pi=B_0a_1B_1a_2B_2a_3B_3\cdots$,
again relaxing the condition that its last element should be smaller
than $a_1$. This gives the upper bound
$e^z\exp(S(z)+e^z-1)=\exp(S(z)+e^z+z-1)$.

\ms

Now we use a similar reasoning to obtain a lower bound. We have seen
that $b_k(z)$ counts blocks of the form $a_iw_{i,1}w_{i,2}\cdots
w_{i,k}$, where each $w_{i,j}$ is a decreasing word starting above
$a_i$ and ending below it. If $k\ge1$, using the notation
$w_{i,k}=w^+_{i,k}w^-_{i,k}$ to separate the elements that are
bigger than $a_i$ from those that are smaller, we can move the last
part of the block to the beginning and write
$L_i:=w^-_{i,k}a_iw_{i,1}w_{i,2}\cdots w^+_{i,k}$. Similarly, a
block of the form $a_iw_{i,0}w_{i,1}w_{i,2}\cdots w_{i,k}$ like the
ones counted by $c_k(z)$ with $k\ge1$ can be reordered as
$L'_i:=w^-_{i,k}a_iw_{i,0}w_{i,1}w_{i,2}\cdots w^+_{i,k}$. The EGF
that counts sets of pieces of the forms given by $L_i$ and $L'_i$ is
$$\exp\left(\sum_{k\ge1}b_k(z)+\sum_{k\ge1}c_k(z)\right)=\exp(S(z)).$$
Ordering the pieces of such a set by decreasing order of the
$a_i$, the sequence that they form by juxtaposition is a
$12\mn34$-avoiding permutation. Besides, no such permutation is
obtained in more than one way by this construction. However,
notice that not every $12\mn34$-avoiding permutation is produced
by this process, hence this construction gives only a lower bound.
\end{proof}

\ms

The decomposition of $12\mn34$-avoiding permutations given in the
proof of Proposition~\ref{prop:12_34} can be generalized to
permutations avoiding a pattern of the form $12\mn\sigma$. If
$\sigma=\sigma_1\sigma_2\cdots\sigma_k\in\S_k$ is a consecutive
pattern, $12\mn\sigma$ denotes the generalized pattern
$12\mn(\sigma_1+2)(\sigma_2+2)\cdots(\sigma_k+2)$.

Any permutation $\pi$ that avoids $12\mn\sigma$ can be uniquely
decomposed as $\pi=B_0a_1B_1a_2B_2a_3B_3\cdots$, where $a_1$ and the
element preceding it form the first ascent of $\pi$, $a_2$ and the
element preceding it form the first ascent such that $a_2<a_1$,
$a_3$ and the element preceding it form the first ascent such that
$a_3<a_2$, and so on. Then, by definition, $B_0$ is a non-empty
decreasing word whose last element is less than $a_1$, and each
$B_i$ with $i\ge1$ can be written uniquely as a sequence
$B_i=w_{i,0}U_{i,1}w_{i,1}U_{i,2}w_{i,2}\cdots U_{i,r_i}w_{i,r_i}$
for some $r_i\ge1$ ($r_i$ can be 0 if $w_{i,0}$ is nonempty) with
the following properties:\ben
\item each $w_{i,j}$ is a decreasing word all of whose elements are less than $a_i$,
\item each $U_{i,j}$ is a nonempty permutation avoiding $\sigma$,
all of whose elements are greater than $a_i$,
\item $w_{i,j}$ is nonempty for $j\ge1$,
\item the last element of $B_i$ is less than $a_{i+1}$.
\een From this decomposition the following result follows
immediately.

\begin{prop} If $\sigma,\tau$ are two consecutive patterns
satisfying $A_\sigma(z)=A_\tau(z)$, then
$A_{12\mn\sigma}(z)=A_{12\mn\tau}(z)$.
\end{prop}

The structure of $21\mn\sigma$-avoiding permutations (defined
analogously) can be described using the same ideas, and it is not
hard to see that the following result holds as well.

\begin{prop}
If $\sigma$ is a consecutive pattern, then
$A_{12\mn\sigma}(z)=A_{21\mn\sigma}(z)$.
\end{prop}

\section{The pattern $1\mn23\mn4$}\label{sec:1_23_4}

Similarly to what we did for the pattern $12\mn34$, analyzing the
structure of permutations avoiding $1\mn23\mn4$ we can give lower
and upper bounds for the numbers $\alpha_n(1\mn23\mn4)$. Let
$\Ce(z):=\sum_{n\ge0}\C_n\frac{z^n}{n!}$ be the EGF for the Catalan
numbers.

\begin{prop}\label{prop:1_23_4}
We have that $$\frac{1}{2}\int_0^ze^{2e^y-2}\,dy-\frac{z}{2} <
A_{1\mn23\mn4}(z)< \Ce(e^z-1).$$
\end{prop}

\begin{figure}[hbt] \fig{file=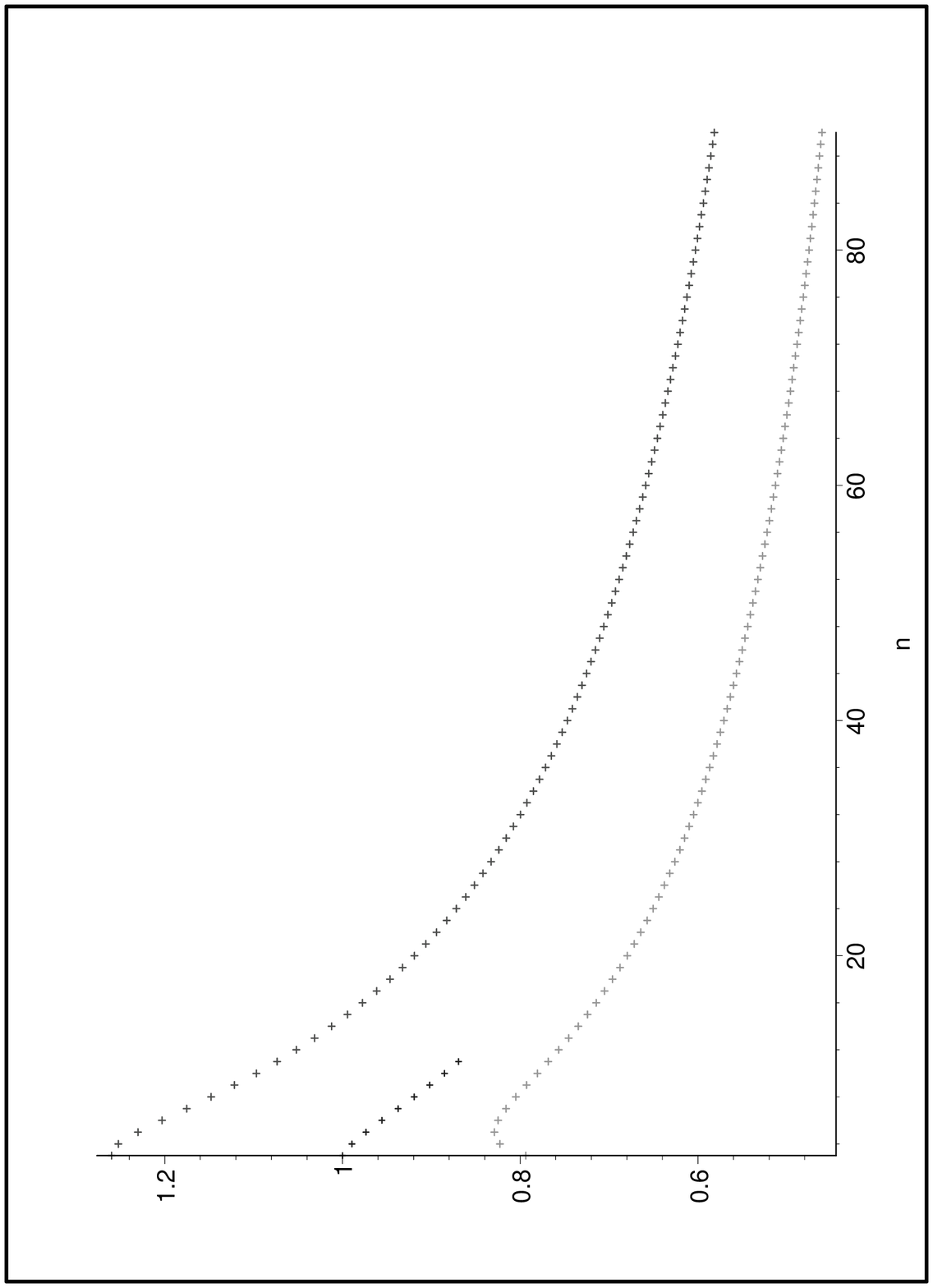,height=12.5cm,angle=-90}
\caption{\label{fig:1bounds4} The first values of
$\sqrt[n]{\alpha_n(1\mn23\mn4)/n!}$ between the lower and the
upper bound given by Proposition~\ref{prop:1_23_4}.}
\end{figure}

Writing $\frac{1}{2}\int_0^ze^{2e^y-2}\,dy-\frac{z}{2}=\sum l_n
\frac{z^n}{n!}$ and $\Ce(e^z-1)=\sum u_n \frac{z^n}{n!}$ to denote
the coefficients of the series giving the lower and the upper
bound respectively, then the values of $\sqrt[n]{l_n/n!}$ and
$\sqrt[n]{u_n/n!}$ for $n\le90$ are plotted in
Figure~\ref{fig:1bounds4}, bounding the values of
$\sqrt[n]{\alpha_n(1\mn23\mn4)/n!}$ for $n\le11$.

Note that the lower bound implies that
$\alpha_n(1\mn23\mn4)\gg\B_n$, since $e^{2e^z-2}\gg e^{e^z-1}$.

\begin{proof}
Let $\pi$ be a permutation that avoids $1\mn23\mn4$. Let
$a_1>a_2>a_3>\cdots>a_r$ be the left-to-right minima of $\pi$, and
let $b_1>b_2>b_3>\cdots>b_s$ be its right-to-left maxima (recall
that $\pi_i$ is a right-to-left maximum of $\pi$ if $\pi_j<\pi_i$
for all $j>i$). Then, marking the positions of the left-to-right
minima and right-to-left maxima, we can write
$\pi=c_1w_1c_2w_2\cdots c_{r+s-1}w_{r+s-1}c_{r+s}$, where
$c_i\in\{a_1,a_2,\ldots,a_r,b_1,b_2,\ldots,b_s\}$ for all $i$ (in
fact the number of $c_i$'s could be less than $r+s$ if some element
is simultaneously a left-to-right minimum and a right-to-left
maximum). Note that $c_1=a_1$ and $c_{r+s}=b_s$. Now, the condition
that $\pi$ avoids $1\mn23\mn4$ is equivalent to the fact that each
$w_i$ is a (possibly empty) decreasing word. Indeed, if there was an
ascent inside one of the $w_i$, then together with the closest
left-to-right minimum to the left of $w_i$ and the closest
right-to-left maximum to the right of $w_i$, it would form an
occurrence of $1\mn23\mn4$. On the other hand, it is clear that if
all $w_i$ are decreasing, then no such occurrence can exist.

We use this decomposition to obtain upper and lower bounds for
$\alpha_n(1\mn23\mn4)$. Let us first show the lower bound. For that
we count only a special type of $1\mn23\mn4$-avoiding permutations,
namely the ones where all the left-to-right minima come before all
the right-to-left maxima. Such a $\pi$ can be written as
$\pi=a_1w_1a_2w_2\cdots a_rw_r b_1w_{r+1}b_2w_{r+2}\cdots
w_{r+s-1}b_s$, where for $1\le i\le r$ the elements of the
decreasing words $w_i$ have values between $a_i$ and $b_1$, and for
$r\le i\le r+s-1$ the elements of $w_i$ have values between $a_r$
and $b_{i+1}$. The EGF for the part $a_1w_1a_2w_2\cdots
a_{r-1}w_{r-1}$ is $e^{e^z-1}$, since it is an arbitrary
$1\mn23$-avoiding permutation (see the example following
Proposition~\ref{prop:1sigma}). Similarly, the EGF for the part
$w_{r+1}b_2w_{r+2}\cdots w_{r+s-1}b_s$ is also $e^{e^z-1}$ (it can
be viewed as a set of blocks of the form $w_{r+i}b_{i+1}$, each one
contributing $e^z-1$, arranged by decreasing order of the $b_i$'s).
The decreasing word $w_r$ contributes $e^z$. Now, to get the EGF for
the whole permutation $a_1w_1a_2w_2\cdots a_rw_r
b_1w_{r+1}b_2w_{r+2}\cdots w_{r+s-1}b_s$ we use the boxed product
construction to require that the biggest element of the block is
$b_1$ and the smallest one is $a_r$. The EGF that we obtain is
$$\int_0^z\int_0^y e^{e^t-1} \left(\frac{d}{dt}t\right) e^t
\left(\frac{d}{dt}t\right)e^{e^t-1}\,dt\,dy=\frac{1}{2}\int_0^z
(e^{2e^y-2}-1)\,dy,$$ which gives a lower bound for the
coefficients of $A_{1\mn23\mn4}(z)$.

To find the upper bound, consider first permutations of the form
$\pi=c_1w_1c_2w_2\cdots c_{r+s-1}w_{r+s-1}c_{r+s}$ where all the
$w_i$ are empty. Such permutations, where every element is either
a left-to-right minimum or a right-to-left maximum, are precisely
those avoiding $1\mn2\mn3$, which are counted by the Catalan
numbers. Thus, the EGF for such permutations is $\Ce(z)$.

The next step is to insert a decreasing word $w_i$ after each $c_i$.
If $c_i$ is a left-to-right minimum, we require that the elements of
$w_i$ are bigger than $c_i$, so the EGF for the block $c_iw_i$ is
$e^z-1$. We omit the requirement that the elements of $w_i$ have to
be smaller than the nearest right-to-left maximum to the right of
$w_i$; this is why we only get an upper bound. Similarly, if $c_j$
is a right-to-left maximum, we require that the elements of $w_j$
are smaller than $c_j$, so the EGF for the block $c_jw_j$ is also
$e^z-1$. We also omit the requirement that after the last
right-to-left maximum there is no decreasing word. Replacing each
$c_i$ for a block $c_iw_i$ as just described translates in terms of
generating functions into substituting $e^z-1$ for the variable $z$
in $\Ce(z)$. This gives the stated upper bound.
\end{proof}

The upper bound given in the above proposition yields the
following corollary.

\begin{corollary} We have that
$$\lim_{n\rightarrow\infty}\left(\frac{\alpha_n(1\mn23\mn4)}{n!}\right)^{1/n}=0.$$
\end{corollary}

\begin{proof}
The power series $\Ce(z)$ can be bounded by
$$\Ce(z)<\sum_{n\ge0}4^n\frac{z^n}{n!}=e^{4z},$$
which converges for all $z$. Therefore, so does $\Ce(e^z-1)$, which
is an upper bound for $A_{1\mn23\mn4}(z)$. The result follows now
from the fact that if $\sum_n f_n z^n$ is an analytic function in
the whole complex plane, then
$\lim_{n\rightarrow\infty}\sqrt[n]{f_n}=0$ (see \cite[Chapter
4]{FlSe98} for a discussion).
\end{proof}

\ms

If $\sigma=\sigma_1\sigma_2\cdots\sigma_{k-2}\in\S_{k-2}$ is a
consecutive pattern, let $1\mn\sigma\mn k$ denote the generalized
pattern $1\mn(\sigma_1+1)(\sigma_2+1)\cdots(\sigma_{k-2}+1)\mn k$.
The decomposition of $1\mn23\mn4$-avoiding permutations given in the
proof of the above proposition can be generalized to permutations
avoiding any pattern of the form $1\mn\sigma\mn k$.

Any permutation $\pi$ that avoids $1\mn\sigma\mn k$ can be uniquely
decomposed as $\pi=c_1w_1c_2w_2\cdots c_{m-1}w_{m-1}c_m$, where the
$c_i$ are all the left-to-right minima and right-to-left maxima of
$\pi$, and each $w_i$ is a permutation that avoids $\sigma$, all of
whose elements are bigger than the closest left-to-right minimum to
its left and smaller than the closest right-to-left maximum to its
right.

Using exactly the same reasoning as in the proof of
Proposition~\ref{prop:1_23_4}, we obtain the following lower and
upper bounds for the numbers $\alpha_n(1\mn\sigma\mn k)$.

\begin{prop}\label{prop:1sigmak} Let $\sigma\in\S_{k-2}$ be a consecutive pattern, and
let $1\mn\sigma\mn k$ be defined as above. Then, $$ \int_0^z\int_0^u
e^{2\int_0^y A_\sigma(t)dt+y}\,dy\,du < A_{1\mn\sigma\mn k}(z) <
\Ce\left(\int_0^z A_\sigma(t)\,dt\right).$$
\end{prop}

\begin{corollary} With the same definitions as in the above
proposition,
$$\lim_{n\rightarrow\infty}\left(\frac{\alpha_n(1\mn\sigma\mn k)}{n!}\right)^{1/n}=
\lim_{n\rightarrow\infty}\left(\frac{\alpha_n(\sigma)}{n!}\right)^{1/n}.$$
\end{corollary}

\begin{proof}
The upper and lower bounds for $A_{1\mn\sigma\mn k}(z)$ given in
Proposition~\ref{prop:1sigmak} are analytic functions of
$A_\sigma(z)$, since essentially they only involve exponentials
and integrals. Therefore, $A_{1\mn\sigma\mn k}(z)$ and
$A_\sigma(z)$ have the same radius of convergence, hence the
limits above coincide.
\end{proof}

Finally, the following proposition is an immediate consequence of
the structure of  $1\mn\sigma\mn k$-avoiding permutations discussed
above. In particular, it implies that
$A_{1\mn23\mn4}(z)=A_{1\mn32\mn4}(z)$.

\begin{prop}
If $\sigma,\tau$ are two consecutive patterns in $\S_{k-2}$
satisfying $A_\sigma(z)=A_\tau(z)$, then $A_{1\mn\sigma\mn
k}(z)=A_{1\mn\tau\mn k}(z)$.
\end{prop}

\section{Other patterns}\label{sec:other}

In Section~\ref{sec:1_23_4} we have proved that
$\alpha_n(1\mn23\mn4)\gg\B_n$ and that $\alpha_n(1\mn23\mn4)\ll c^n
n!$ for any constant $c>0$. For the pattern $12\mn34$, we showed in
Section~\ref{sec:12_34} that the analogue to the first statement
holds as well, and the second one seems to be true from numerical
computations. It remains as an open problem to describe precisely
the asymptotic behavior of $\alpha_n(\sigma)$ for these two
patterns, and for several remaining generalized patterns of length
4.

\begin{figure}[hbt] \fig{file=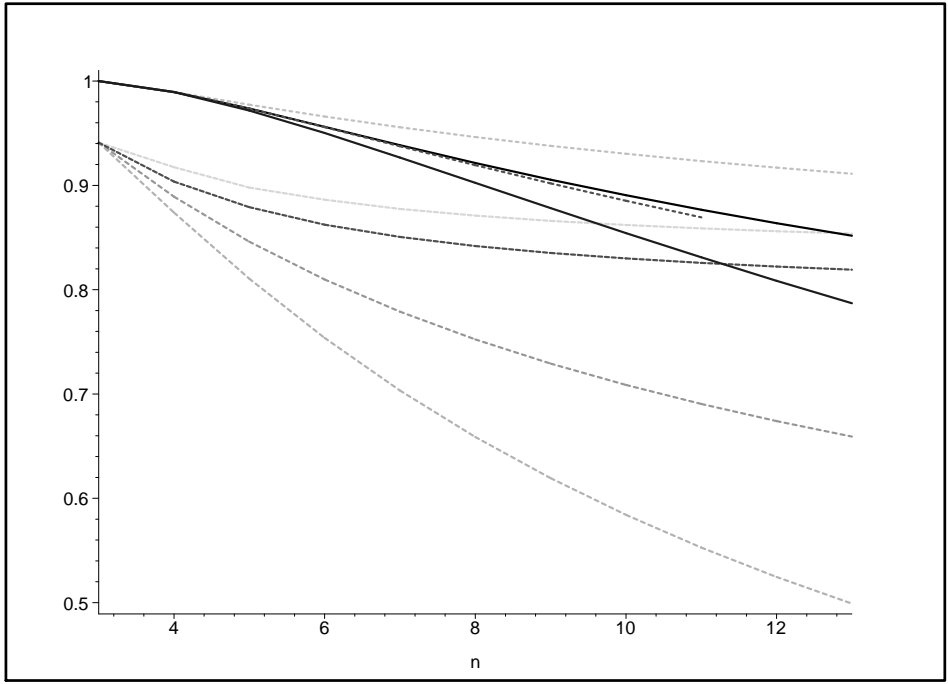,height=12.5cm,angle=-90}
\caption{\label{fig:bigraph} The first values of
$\sqrt[n]{\alpha_n(\sigma)/n!}$ for several generalized patterns
$\sigma$.}
\end{figure}

In Figure~\ref{fig:bigraph} we have plotted the initial values
(connected by lines) of the sequences
$\sqrt[n]{\alpha_n(\sigma)/n!}$ for other cases that appear to have
some interest. The two dotted lines at the bottom of the graph
correspond to the sequences $\sqrt[n]{\C_n/n!}$ and
$\sqrt[n]{\B_n/n!}$, which are known to tend to 0 as $n$ goes to
infinity. The two dashed lines that start at the same point (around
$0.941$) and tend to a constant correspond to the sequences
$\sqrt[n]{\alpha_n(132)/n!}$ and $\sqrt[n]{\alpha_n(123)/n!}$, for
which their limits are known by Theorem~\ref{th:EN} to be
$0.7839769$ and $0.8269933$ respectively. Among the lines starting
at 1, the two dotted ones correspond to the patterns $1\mn23\mn4$
(the lower line) and $12\mn34$ (the upper line) discussed in the
previous sections.

Of the two solid lines, the one below corresponds to the pattern
$3\mn14\mn2$. This pattern has a special interest because all of its
subpatterns of length~3 are among those in part~(ii) of
Proposition~\ref{prop:claesson}. Since it does not contain any of
the patterns in part~(i), we cannot say that
$\alpha_n(3\mn14\mn2)\ge\B_n$ for all $n$. In fact, comparing the
slopes in~Figure~\ref{fig:bigraph} it seems quite plausible that
$\alpha_n(3\mn14\mn2)$ grows more slowly than $\B_n$, and proving
this is an interesting open question. The other solid line in the
plot corresponds to the pattern $13\mn24$, for which we do not know
the asymptotic behavior either.

\ms

This paper is the first attempt to study the asymptotic behavior of
the numbers $\alpha_n(\sigma)$ where $\sigma$ is an arbitrary
generalized pattern. Despite the fact that we have been unable to
provide a precise description of this behavior in most cases, we
hope that our work shows the intricateness of the problem and the
amount of questions that it opens. The main goal of further research
in this direction would be to give a complete classification of all
generalized patterns according to the asymptotic behavior of
$\alpha_n(\sigma)$ as $n$ goes to infinity.

Another interesting open problem is to find the value of
$\lim_{n\rightarrow\infty}\sqrt[n]{\alpha_n(\sigma)/n!}$ for
patterns $\sigma$ in Case 2, for which this limit is known to be a
constant. The analogous problem for patterns in Case 1, namely
finding $\lim_{n\rightarrow\infty}\sqrt[n]{\alpha_n(\sigma)}$ for
classical patterns $\sigma$, is a current direction of research as
it remains open for most patterns as well (see \cite{Bon97,B3}).

\end{document}